\pdfoutput=1

\documentclass{article}
\usepackage[T1]{fontenc}
\usepackage[UKenglish]{babel}

\usepackage{Factory, Math, Theorema, Classic, Styl}

\newcommand{\eqar}{\ar@{=} }

\newcommand{\Zp}{\Z/p}					
\newcommand{\freenum}[2][\Z]{#1\binom{#2}{-}}

\usepackage{fancyhdr}						
\begin{document}

\lefthyphenmin=2 \righthyphenmin=2

\titul{BINOMIAL RINGS: AXIOMATISATION, TRANSFER AND CLASSIFICATION}
\auctor{Qimh Richey Xantcha\thanks{\textsc{Qimh Richey Xantcha}, Stockholm University: \texttt{qimh@math.su.se}}}
\datum{\today}
\maketitle

\bigskip 


\bigskip 

\begin{argument} 
\noindent
Hall's \emph{binomial rings}, rings with binomial coefficients, are given an axiomatisation and proved identical to 
the \emph{numerical rings} studied by Ekedahl. The \emph{Binomial Transfer Principle} is established, 
enabling combinatorial proofs of algebraical identities. The finitely generated binomial rings 
are completely classified. An application to modules over binomial rings is given. 

\end{argument}

\bigskip 

\noindent 
The abstract study of binomial coefficients seems to have been initiated by Hall \cite{Hall} (\S 6),  
who introduced the concept of binomial rings in connexion with his ground-breaking work on nilpotent groups. 
The definition is simple. A \textbf{binomial ring} is a commutative, unital ring $R$ which is torsion-free and 
closed in $\Q\otimes R$ under the ``formation of binomial coefficients'': 
$$ 
r\mapsto \binom{r}{n} = \frac{r(r-1)\cdots (r-n+1)}{n!}.
$$
(Hall required $R$ to be a domain.) 

Binomial rings naturally manifest themselves in the theories of integer-valued polynomials, Witt vectors 
and $\lambda$-rings; 
to name but a few. We refer the reader to Elliott's article \cite{Elliott} and the lucid monograph \cite{Yau} by Yau, 
where these topics have been expounded upon. 

More recently, binomial rings have turned out to 
form the natural frame-work for discussing polynomial maps and functors of modules;
see \cite{PM} and \cite{PF}.

What is known on binomial rings stems principally from a recent paper \cite{Elliott} by Elliott, which, 
 in particular, aims to elucidate the connexion between binomial rings and $\lambda$-rings.  
Let us compile a list of their most important properties.

\begin{enumerate}
\item \label{P: Free} The free binomial ring on the set $X$ is the ring 
$$
\{f\in\Q[X] \mid f(\Z^X)\subseteq \Z \}
$$
of integer-valued polynomials on $X$. 
(\cite{Elliott}, Proposition 2.1) 

\item  \label{P: Fermat} 
The following conditions on a commutative, unital ring $R$ are equivalent (\cite{Elliott}, Theorem 4.1, 4.2):
\begin{enumerate} 
\item $R$ is the quotient of a  binomial ring. 
\item The elements of $R$ satisfy every integer \emph{polynomial congruence} universally true for the integers. 
\item $a(a-1)\cdots(a-n+1) $ is divisible by $n!$ for every  $n\in \N$.
\item \emph{Fermat's Little Theorem} holds: $a^p\equiv a \mod pR$ for every prime $p$. 
\item The Frobenius map $a\mapsto a^p$ is the identity on $R/pR$ for every prime~$p$. 
\item $R/pR$ is reduced for every prime $p$, and each of its residue fields is isomorphic to $\Zp$. 
\end{enumerate} 
Adjoining, in each case, the assumption that $R$ lack torsion, these transform into criteria for $R$ to be binomial. 

\item The binomial property is preserved under the following constructions: localisation, direct products, 
tensor products,  filtered inductive and projective limits. (\cite{Elliott}, Propositions 5.1, 5.4, 5.5)

\item The inclusion functor from binomial rings to rings has both a left and a right adjoint. 
(\cite{Elliott}, Theorems 7.1, 9.1)

\item Binomial rings are equivalent to $\lambda$-rings with trivial Adams operations. 
(\cite{Wilkerson}, Proposition 1.2; \cite{Elliott}, Proposition 8.3)

\item \emph{The Binomial Theorem.} Let $R$ be binomial 
and let $A$ be a commutative algebra over $R$ which is complete with respect to the ideal $I$. 
The equation $$ (1+x)^r = \sum_{n=0}^\infty \binom rn x^n $$ defines an $R$-module structure on the abelian 
group $(1+I,\cdot)$.
(\cite{Elliott}, Proposition 11.1)

\end{enumerate}

Ekedahl \cite{TE}, preferring the axiomatic approach, proposed six axioms intended to capture the properties
of binomial coefficients. 
He appears not to have been familiar with 
the work of Hall and never proved his \emph{numerical rings} (as he called them) 
to be equivalent to Hall's binomial rings. 
Rectifying this is one object of the present paper. 
Indeed, we not only justify, but improve upon Ekedahl's axioms, 
dropping the ghastliest axiom (the sixth): 

\begin{inttheorem}[\ref{S: Coincide}] 
In order for a commutative, unital ring to be binomial, it is necessary and sufficient that it be equipped
with unary operations $r\mapsto \binom{r}{n}$ ($n\in\N$), satisfying the following five axioms.
\begin{Romanlist}
\item $\displaystyle \binom{a+b}{n} = \sum_{p+q=n} \binom{a}{p} \binom{b}{q}$. 
\item $\displaystyle \binom{ab}{n} = \sum_{m=0}^{n} \binom{a}{m} \sum_{\substack{q_1+\cdots+q_m=n \\ q_i\geq 1}} \binom{b}{q_1}\cdots\binom{b}{q_m}$.
\item $\displaystyle \binom{a}{m}\binom{a}{n} =  \sum_{k=0}^{n} \binom{a}{m+k}\binom{m+k}{n}\binom{n}{k}$.
\item $\displaystyle \binom1n=0$\quad when $n\geq 2$. 
\item $\displaystyle \binom a0=1$\quad and\quad $\displaystyle \binom a1=a$.
\end{Romanlist} \hfill
\end{inttheorem}

Ekedahl's sixth axiom was a (non-explicit) formula for reducing an iterated binomial coefficient
$\binom{\binom rm}{n}$ to simple ones. 
Surprisingly, such a formula will turn out to be a consequence of the five axioms listed. 
(No explicit formula for iterated binomial coefficients 
appears to exist. Golomb \cite{Golomb} has summarised what is known on the subject.)

We next prove the following \emph{Transfer Principle}, formally sanctioning  
combinatorial proofs of algebraical identities in binomial rings. Compare property 
\ref{P: Fermat}(b) above.   

\begin{inttheorem}[\ref{S: Transfer} (The Binomial Transfer Principle)] 
A binomial polynomial identity 
universally valid in $\Z$ is valid 
in every binomial ring. \end{inttheorem}

Our most important result  
is the complete classification of all finitely generated binomial rings.

\begin{inttheorem}[\ref{S: Classification} (The Structure Theorem for Finitely Generated Binomial Rings)]
Let $R$ be a finitely generated binomial ring. 
There exist unique positive, simply composite integers $m_1,\dots,m_k$ such that 
$$ 
R\cong \Z[m_1^{-1}]\times \cdots\times \Z[m_k^{-1}]. 
$$ 
\end{inttheorem}

The classification is used in the treatment of torsion-free modules in Section \ref{A: Modules}.

\section{Definition and Examples}

Hall's original definition, as found in \cite{Hall} (Section 6), reads as follows.

\begin{definition} 
A commutative ring $R$ with unity is a \textbf{binomial ring} if it is torsion-free\footnote{The 
word \emph{torsion} will, here and elsewhere, be taken to mean \emph{$\Z$-torsion}.} and 
closed in $\Q\otimes R$ under the operations 
$$ r\mapsto \frac{r(r-1)\cdots (r-n+1)}{n!}.$$ 
\end{definition}

\begin{example}
In any $\Q$-algebra, binomial coefficients may be defined by the usual formula. 
\end{example}

\begin{example} 			\label{Ex: Num Z[m^{-1}]}
For any integer $m\neq 0$, the ring $\Z[m^{-1}]$ is binomial. Since it inherits the binomial coefficients from $\Q$, 
it is simply a question of verifying 
closure under the formation of binomial coefficients. Because
$$
\binom{\frac{a}{f}}{n} = \frac{\frac af(\frac af -1)\cdots (\frac af -(n-1)) }{n!} 
= \frac{a(a-f)\cdots(a-(n-1)f)}{n!f^n},
$$
it will suffice to prove that whenever $p^i\mid n!$, but $p\nmid b$, then 
$$ 
p^i\mid (a+b)(a+2b)\cdots (a+nb).
$$

To this end, let
$$ 
n = c_mp^m + \cdots + c_1 p + c_0 , \qquad 0\leq c_i\leq p-1,
$$
be the base $p$ representation of $n$. For fixed $k$ and $0\leq d<c_k$, the numbers
\begin{equation} \label{Eq: Z[m^-1] A}
a + (c_m p^m + \cdots + c_{k+1}p^{k+1} + dp^k + i)b , \qquad 1\leq i\leq p^k,
\end{equation}
will form a set of representatives for the congruence classes modulo $p^k$, as will of course the numbers 
\begin{equation} \label{Eq: Z[m^-1] B}
c_m p^m + \cdots + c_{k+1}p^{k+1} + dp^k + i , \qquad 1\leq i\leq p^k.
\end{equation}
Note that if $x\equiv y\mod p^k$ and $j\leq k$, then $p^j\mid x$ iff $p^j\mid y$. 
Hence there are at least as many factors $p$ among the numbers \eqref{Eq: Z[m^-1] A} as among the numbers \eqref{Eq: Z[m^-1] B}. The claim now follows. 
\end{example}

\begin{example}
Being given by rational polynomials, the operations $r\mapsto \binom{r}{n}$ give continuous maps 
$\Q_p\to\Q_p$ in the $p$-adic topology. 
It should be well known that $\Z$ is dense in the ring $\Z_p$, and that $\Z_p$ is closed in $\Q_p$. 
Since the binomial coefficients leave $\Z$ invariant, 
the same must be true of $\Z_p$, which is thus binomial. 

This provides an alternative proof of the fact that $\Z[m^{-1}]$ is closed under binomial coefficients. 
For this is evidently true of the localisations 
$$
\Z_{(p)}=\Q\cap\Z_p,
$$ 
and therefore also for 
\[
\Z[m^{-1}] = \bigcap_{p\nmid m} \Z_{(p)}. \qedhere
\] 
\end{example}

\section{Axioms}

We next present, with minor modifications, Ekedahl's axioms for numerical rings (\cite{TE}, Definition 4.1), with the 
notable exception of the sixth. They have been amended so as to yield explicit formul\ae. 

\begin{definition}
A \textbf{numerical ring} is a commutative ring with unity equipped with unary operations 
$$
r\mapsto \binom{r}{n}, \qquad n\in\N;
$$ 
called \textbf{binomial coefficients}, subject to the following axioms.
\begin{Romanlist}
\item $\displaystyle \binom{a+b}{n} = \sum_{p+q=n} \binom{a}{p} \binom{b}{q}$. 
\item $\displaystyle \binom{ab}{n} = \sum_{m=0}^{n} \binom{a}{m} \sum_{\substack{q_1+\cdots+q_m=n \\ q_i\geq 1}} \binom{b}{q_1}\cdots\binom{b}{q_m}$.
\item $\displaystyle \binom{a}{m}\binom{a}{n} =  \sum_{k=0}^{n} \binom{a}{m+k}\binom{m+k}{n}\binom{n}{k}$.
\item $\displaystyle \binom1n=0$\quad when $n\geq 2$. 
\item $\displaystyle \binom a0=1$\quad and\quad $\displaystyle \binom a1=a$.
\end{Romanlist} 
\end{definition}

For $a,b\in\Z$, these identities embody solutions to certain problems of enumerative combinatorics, 
which we leave for the reader to formulate.

It follows easily from Axioms \textsc{i}, \textsc{iv} and \textsc{v} that, 
when the functions $\binom -n$ are evaluated on natural multiples of unity, we 
retrieve the ordinary binomial coefficients, namely 
$$
\binom{m\cdot 1}{n}= \frac{m(m-1)\cdots (m-n+1)}{n!}\cdot 1, \quad m\in\N.
$$ 

Since $\binom {n\cdot 1}n=1$, but $\binom0n=0$ unless $n=0$, a numerical ring has necessarily characteristic $0$, 
and so contains $\Z$. 
No confusion arises from writing $\binom{r}{n}$ when $r\in\Z\subseteq R$:

\begin{lemma} 				\label{L: Elementary} 
The following formul\ae{} are valid in a numerical ring: 
\begin{palphalist}
\item 			\label{L: Elementary a}
$\displaystyle \binom{r}{n}=\frac{r(r-1)\cdots (r-n+1)}{n!}$ \quad when $r\in \Z$. 
\item 			\label{L: Elementary b}
$\displaystyle n!\binom{r}{n}= r(r-1)\cdots (r-n+1)$ \quad for any $r$.
\item 			\label{L: Elementary c}
$\displaystyle n\binom{r}{n}=(r-n+1)\binom{r}{n-1}$ \quad for any $r$.
\end{palphalist}
\end{lemma}

\begin{proof} The map 
$$ 
\phi\colon (R,+) \to (1+tR[[t]],\cdot), \qquad r\mapsto \sum_{n=0}^\infty \binom rn t^n ,
$$
is, by Axioms \textsc{i} and \textsc{v}, a group homomorphism. Therefore, when $r\in \Z$, 
$$ 
\phi(r) = \phi(1)^r = (1+t)^r, 
$$
which expands as usual (with ordinary binomial coefficients) by the Binomial Theorem. This proves part (a). 
(An inductive proof would also work.)

To prove parts (b) and (c), we proceed differently. By Axiom \textsc{iii},
\begin{align*}
r\binom{r}{n-1} &=\binom{r}{n-1}\binom{r}{1} = \sum_{k=0}^1 \binom{r}{n-1+k}\binom{n-1+k}{1}\binom{1}{k} \\
& = \binom{r}{n-1}\binom{n-1}{1}\binom{1}{0} + \binom{r}{n}\binom{n}{1}\binom{1}{1} \\ &= (n-1)\binom{r}{n-1} + n\binom{r}{n},
\end{align*}
which reduces to (c). Part (b) then follows inductively from (c). 
\end{proof}

Our present objective will be showing that numerical and binomial rings coincide. 
It follows that the numerical structure on a given ring is always unique. 

\begin{lemma} 		\label{L: Num. tfree A}
Let $m$ be an integer. If $p$ is prime and $p^l\mid m$, but $p\nmid k$, then $p^l\mid \binom{m}{k}$. 
\end{lemma}

\begin{proof}
The number $p^l$ divides the right-hand side of 
$$
k\binom{m}{k}=m\binom{m-1}{k-1},
$$ and therefore also the left-hand side. But $p^l$ is relatively prime to $k$, 
so in fact $p^l\mid \binom{m}{k}$.
\end{proof}

\begin{lemma} 		\label{L: Num. tfree B}
Let $m_1,\dots,m_k$ be natural numbers. 
If 
$$
m_1 + 2m_2 + 3m_3 + \dots + km_k
$$ 
is prime, then 
$$ 
m_1+\dots+m_k \mid \binom{m_1+\dots+m_k}{m_1,\dots,m_k}, 
$$ 
unless $m_2=\dots=m_k=0$. 
\end{lemma}

\begin{proof} 
Consider a prime power $p^l\mid m_1+\dots+m_k$. 
If $p\nmid m_j$, then 
$$ 
\binom{m_1+\dots+m_k}{m_1,\dots,m_k} = \binom{m_1+\dots+m_k}{m_j} \binom{\sum_{i\neq j}m_i}{\{m_i\}_{i\neq j}} 
$$
is divisible by $p^l$ according to Lemma \ref{L: Num. tfree A}, and the assertion follows. 
Supposing all $m_i$ to be divisible by $p$, it follows that 
$$
p\mid m_1 + 2m_2 + 3m_3 + \dots + km_k.
$$ 
Since this latter number is itself prime, $m_1=p$ and $m_2=\dots=m_k=0$.
\end{proof}

\begin{lemma} 			\label{L: Num. tfree C}
Let $R$ be a numerical ring. Let $r\in R$ and $m,n\in\N$. If $nr=0$, then also $mn\binom rm=0$. 
\end{lemma}

\begin{proof} 
This follows inductively, since  
$$ 
mn\binom{r}{m} = n(r-m+1)\binom{r}{m-1} = -n(m-1)\binom{r}{m-1},
$$ 
using part \ref{L: Elementary c} of Lemma \ref{L: Elementary}.
\end{proof}

\begin{theorem} 
Numerical rings are torsion-free. 
\end{theorem}

\begin{proof} 
Let $r\in R$ (numerical) and let $n$ be a prime number.  
Supposing $nr=0$, we calculate
$$
0=\binom 0n = 
\binom{nr}{n} = 
\sum_{m=0}^{n} \binom{r}{m} \sum_{\substack{q_1+\cdots+q_m=n \\ q_i\geq 1}} \binom{n}{q_1}\cdots\binom{n}{q_m}.
$$
Assigning $m_j$, for $1\leq j\leq n$, to stand for the number of $q_i$ that are equal to 
$j$ will transform the sum into
\begin{equation}			\label{E: Torsion-free}
0 = 
\sum_{m=0}^{n} \binom{r}{m} \sum_{\substack{m_1+\cdots+m_n=m \\ m_1+2m_2+\cdots+ nm_n = n}} 
\binom{m}{m_1,\dots,m_n} \binom{n}{1}^{m_1}\binom{n}{2}^{m_2}\cdots\binom{n}{n}^{m_n}.
\end{equation}
(Given that the value $j$ occurs exactly $m_j$ times among the $q_i$, values may distributed to 
the variables $q_i$ in $\binom{m}{m_1,\dots,m_n}$ ways, accounting for the multinomial 
coefficient above.)

The inner sum is empty when $m=0$. For $m=1$, it must be that $m_1=\dots=m_{n-1}=0$ and $m_n=1$, and 
the inner sum reduces to $\binom{r}{1}$. 

In so far as $2\leq m\leq n-1$, each term of the inner sum is divisible by $mn$. 
Indeed, $m_j>0$ for some $j<n$, for which
Lemma \ref{L: Num. tfree A} asserts that $n\mid\binom{n}{j}^{m_j}$.
Also, $m_k>0$ for some $k>1$, and then $m\mid\binom{m}{m_1,\dots,m_n}$ by Lemma \ref{L: Num. tfree B}.

In the case $m=n$, obviously $m_2=\cdots=m_n=0$ and $m_1=n$. 
The inner sum will equal $\binom{n}{1}^n$, which is again divisible by $n^2=mn$. 

The terms of \eqref{E: Torsion-free} corresponding to $m\geq 2$ will therefore vanish by  Lemma~\ref{L: Num. tfree C}, 
leaving only the summand $m=1$. But this term is 
simply $\binom r1=r$, which is then equal to zero, and $R$ is torsion-free. 
\end{proof}


\begin{theorem} 			\label{S: Coincide} 
Numerical and binomial rings coincide. 
\end{theorem}

\begin{proof} 
Clearly, binomial rings satisfy the numerical axioms. 
Conversely, the binomial coefficients of a numerical ring fulfil 
$$
n!\binom{r}{n}= r(r-1)\cdots (r-n+1)
$$ 
by part (b) of Lemma \ref{L: Elementary}, which, by the absence of torsion, implies
\[
\binom{r}{n}= \frac{r(r-1)\cdots (r-n+1)}{n!}.		\qedhere
\]
\end{proof}

The appellations \emph{numerical} and \emph{binomial} may thus be treated synonymously.


Let us now resolve the mystery of the missing sixth axiom. 
In $\Z$, there ``exists'' a formula for iterated binomial coefficients: 
\begin{equation} \label{Eq: It BC} \binom{\binom{r}{m}}{n}=\sum_{k=1}^{mn} g_k\binom{r}{k}, \end{equation}
in the sense that there are unique integers $g_k$ making the formula valid for every $r\in\Z$. 
Golomb has examined these iterates in some detail, 
and his paper \cite{Golomb} is brought to an end with the discouraging conclusion: 
\begin{quote}\small No simple reduction formulas have yet been found for the most general case of $\binom{\binom nb}a$. \end{quote}

Note, however, that \eqref{Eq: It BC} is a polynomial identity with rational coefficients, 
by which it must hold in every $\Q$-algebra, 
and therefore in every binomial ring. This proves the redundancy of Ekedahl's original sixth axiom:

\begin{theorem} 			\label{S: Iterated}
The formula 
$$
\binom{\binom{r}{m}}{n}=\sum_{k=1}^{mn} g_k\binom{r}{k}
$$ 
for iterated binomial coefficients is valid in every binomial ring. 
\end{theorem}

\section{Transfer}

%
Let $X$ be a set, and let $E(X)$ be the \emph{term algebra}%
\footnote{The denomination
\emph{term algebra} is borrowed from universal 
algebra; confer Definition II.10.4 of \cite{Universal}.} 
based on $X$. It consists of all finite words that can be formed from the alphabet 
$$ 
X \cup\Set{ +,  -,  \cdot,  0,  1,  \binom-n }{  n\in\N }, 
$$
with binary operations $+$ and $\cdot$, unary operations $-$ and $\binom -n$, 
and nullary operations $0$ and $1$ (constants). 

\begin{definition} 
The ring $\freenum{X}$ is the result of imposing  upon the term algebra the axioms of a 
commutative ring with unity, along with the numerical axioms. 
\end{definition}

The numerical axioms, together with the formula for iterated binomial coefficients 
(Theorem \ref{S: Iterated}) will 
reduce any element of $\freenum{X}$ to a \emph{binomial polynomial} of the form
$$
\sum_{n_1,\dots,n_k=0}^{m} c_{n_1,\ldots,n_k}\binom{x_1}{n_1}\cdots\binom{x_k}{n_k}, 
\qquad x_i\in X, \ c_{n_1,\ldots,n_k}\in \Z. 
$$
Conversely, as is well known, any integer-valued polynomial $f(x)\in\Q[X]$ is given by a unique binomial polynomial, 
and we have provided another description of the free binomial ring on $X$; 
confer property \ref{P: Free} in the introductory section:

\begin{theorem}  
There is an isomorphism 
$$
\freenum{X} \cong \set{f\in\Q[X] }{ f(\Z^X)\subseteq \Z }
$$
with the free binomial ring on the set $X$. 
\end{theorem}

We derive the following very useful corollary. 

\begin{theorem}[The Binomial Transfer Principle]			\label{S: Transfer} 
A binomial polynomial identity universally valid in $\Z$ is valid 
in every binomial ring. 
\end{theorem}

\begin{proof} 
A binomial polynomial $p(x_1,\dots,x_k)$ 
evaluating to zero on all integers is the zero integer-valued polynomial and therefore, 
by the previous theorem, the zero binomial polynomial. 
\end{proof}

\section{Binomial Ideals and Factor Rings}

We now make a short survey of binomial ideals and the associated factor rings. 

\begin{theorem} 
Let $I$ be an ideal of the binomial ring $R$. The equation 
$$ 
\binom{r+I}{n} = \binom{r}{n} + I 
$$
will yield a binomial structure on $R/I$ if and only if $\binom{e}{n}\in I$ 
for every $e\in I$ and $n>0$. 
\end{theorem}

\begin{proof} 
The condition is clearly necessary, as we must have $\binom{e}{n}\equiv \binom{0}{n} \mod I$, for any $e\in I$,
in order for the binomial structure to be well defined. Conversely, assuming 
$\binom{e}{n}\in I$ for every $e\in I$ and $n>0$, the structure is indeed well defined: 
$$ 
\binom{r+e}{n} = \sum_{p+q=n} \binom rp \binom eq \equiv \binom rn \binom e0 = \binom rn \mod I.
$$ 
The numerical axioms in $R/I$ follow immediately from those in $R$. 
\end{proof}

\begin{definition} 
An ideal of a binomial ring satisfying the condition of the previous theorem will be called a \textbf{binomial ideal}. 
\end{definition}

\begin{example}
$\Z$ does not possess any non-trivial binomial ideals, because all its non-trivial factor rings have torsion. 
Neither do the rings $\Z[m^{-1}]$. 
\end{example}

The next theorem provides a kind of converse. 

\begin{theorem} 
Let $R$ be a (commutative, unital) ring with an ideal $I$. 
Suppose $I$ is a vector space over $\Q$ and that $R/I$ is binomial. 
Then $R$ itself is binomial, and $I$ is a binomial ideal. 
\end{theorem}

\begin{proof} 
Since $I$ and $R/I$ are both torsion-free, so is $R$, and there is a commutative diagram with exact rows: 
$$ 
\xymatrix{
0 \ar[r] & I \eqar[d] \ar[r] & R \ar[d] \ar[r] & R/I \ar[d] \ar[r] & 0 \\
0 \ar[r] & \Q\otimes I \ar[r] & \Q\otimes R \ar[r] & \Q\otimes R/I \ar[r] & 0
} 
$$
Also, $\Q\otimes I=I$, and we may identify 
$$
\Q\otimes R/I \cong (\Q\otimes R)/I.
$$
For $r\in R$, we note that
$$
\binom{r+I}{n} = \frac{r(r-1)\cdots (r-n+1)}{n!} +I 
$$
does in fact lie in $R/I$, being binomial. It must then be that 
$$
\frac{r(r-1)\cdots (r-n+1)}{n!}\in R.
$$ 
which is thus binomial, having $I$ as a binomial ideal.
\end{proof}


\section{Finitely Generated Binomial Rings}

We shall now classify the finitely generated binomial rings. 

\begin{lemma} 			\label{L: Ring of fractions}
If a (commutative, unital) 
ring $R$ is torsion-free and finitely generated as an abelian group, its ring of fractions is $\Q\otimes R$. 
\end{lemma}

\begin{proof} 
By the Structure Theorem for Finitely Generated Abelian Groups, $R$ is isomorphic to some $\Z^n$ as an abelian group. 
Let $a\in \Z^n$. Multiplication by $a$ is 
a linear transformation on $\Z^n$, and so may be represented by an integer matrix $A$. 
In so far as $a$ not be a zero-divisor, $A$ will be 
non-singular. It will then have an inverse $A^{-1}$ with rational entries. The inverse of $a$ is given by 
$$ 
a^{-1}=A^{-1}1 \in \Q^n=\Q\otimes R, 
$$ 
where the multiplicative identity $1\in R$ has been written as a column vector. 
\end{proof}

The next lemma is considered well known. 

\begin{lemma} 			\label{L: Algebraic integers}
Let $A$ denote the algebraic integers in the field $K\supseteq \Q$. If $K$ is finitely generated over $\Q$, then 
$A$ is finitely generated over $\Z$. 
\end{lemma}

The subsequent (in)equality of Krull dimensions is supposedly familiar to scholars 
in the fields of commutative algebra or 
algebraic geometry. We are grateful to Prof.~T. Ekedahl for furnishing the proof. 

\begin{theorem} 			\label{S: Dimension} 
Let $R$ be a finitely generated, commutative and unital non-zero ring. Then 
$$ 
\dim \Q \otimes R = \dim R/pR \leq \dim R-1 
$$ 
holds for all but finitely many prime numbers $p$. 

When $R$ is an integral domain of characteristic $0$, 
there is in fact equality for all but finitely many primes $p$. 
\end{theorem}

\begin{proof} 
In the case of positive characteristic $n$, the inequality will hold trivially, for then 
$$ 
\Q\otimes R = 0 = R/pR , 
$$ 
except when $p\mid n$. 

Consider now the case when $R$ is an integral domain of characteristic 0. 
There is an embedding $\phi\colon\Z\to R$ and a corresponding 
dominant morphism 
$$
\Spec\phi\colon\Spec R\to\Spec \Z
$$ 
of integral schemes, which is of finite type. Letting $\Frac D$ denote the field of fractions of the integral domain 
$D$, we may define
\begin{align*}
C_n & = \set{P\in \Spec \Z }{ \dim(\Spec\phi)^{-1}(P)=n } \\
&= \set{P\in\Spec \Z }{ \dim R\otimes\Frac R/P=n } \\
&= \set{(p) }{ \dim R/pR  = n } \cup \{(0) \mid \dim R\otimes \Q = n \}  .
\end{align*}
This latter set, by Chevalley's Constructibility Theorem%
\footnote{This proposition appears to belong to the folklore of 
Algebraic Geometry. An explicit reference is Théorème 2.3 of \cite{Jou}.}%
, will contain a dense, open set in $\Spec\Z$ if $n=\dim R-\dim \Z$. 
Such a set must contain $(0)$ and $(p)$ for all but finitely many primes $p$. For these primes,
$$ 
\dim \Q\otimes R = \dim R/pR = n = \dim R - 1.
$$

Now let $R$ be an arbitrary ring of characteristic 0. 
For any prime ideal $Q$, $R/Q$ will be an integral domain (but not necessarily of characteristic $0$), 
and so we can apply the preceding to obtain
$$ 
\dim\Q\otimes R/Q = \dim R/(Q+pR) \leq \dim R/Q -1, 
$$
for all but finitely many primes $p$. 
The prime ideals of $\Q\otimes R$ are of the form $\Q\otimes Q$, where $Q$ is a prime ideal in $R$. Moreover, 
$$ 
(\Q\otimes R)\big/(\Q\otimes Q)=\Q\otimes R/Q.
$$ 
It follows that  
\begin{align*}
\dim \Q\otimes R &= \max_{Q\in\Spec R} \dim (\Q\otimes R)\big/(\Q\otimes Q) \\
&= \max_{Q\in\Spec R} \dim \Q\otimes R/Q \\
&= \max_{Q\in\Spec R} \dim R/(Q+pR) \\
&= \max_{\overline Q\in \Spec R/pR} \dim (R/pR)\big/\overline Q = \dim R/pR
\end{align*}
for all but finitely many $p$ (the maxima are taken over the finitely many \emph{minimal} prime ideals only). 
In a similar fashion, 
\begin{align*}
\dim \Q\otimes R &= \max_{Q\in\Spec R} \dim (\Q\otimes R)\big/(\Q\otimes Q) \\
&= \max_{Q\in\Spec R} \dim \Q\otimes R/Q \\
&\leq \max_{Q\in\Spec R} \dim  R/Q -1 = \dim R - 1. \qedhere
\end{align*}
\end{proof}

The following Classification Theorem, together with its proof, was communicated to us by Prof.~Ekedahl.

\begin{theorem}[The Structure Theorem for Finitely Generated Binomial Rings] 	\label{S: Classification}
Let $R$ be a finitely generated binomial ring. 
There exist unique positive, square-free 
integers $m_1,\dots,m_k$ such that 
$$ 
R\cong \Z[m_1^{-1}]\times \cdots\times \Z[m_k^{-1}]. 
$$ 
\end{theorem}

\begin{proof} \emph{Case A: $R$ is finitely generated as an abelian group.} 
If $r^n=0$, then, because of Fermat's Little Theorem (property \ref{P: Fermat} in the introduction), 
$r$ is divisible by $p$ for all primes $p>n$. But in $\Z^n$ this can only be if $r=0$; hence $R$ is reduced. 
By Lemma \ref{L: Ring of fractions}, the ring of fractions is $\Q\otimes R$. 
As this is reduced and artinian, being finite-dimensional over $\Q$, 
it splits into a product of fields of characteristic $0$. 

\emph{Case A1: The ring of fractions is a field.} 
In the special case when the fraction ring $\Q\otimes R$ is a field, 
whose ring of algebraic integers we denote by $A$, we examine the subgroup $A\cap R$ of $A$. 
Since $A\subseteq \Q\otimes R$, an arbitrary 
element of $A$ will have an integer multiple lying in $R$. This means $A/(A\cap R)$ is a torsion group. 
Also, the fraction ring $\Q\otimes R$ is finitely generated over $\Q$, so, 
from Lemma \ref{L: Algebraic integers}, we deduce that 
$A$ is finitely generated over $\Z$. 
Because the factor group $A/(A\cap R)$ is both finitely generated and torsion, it is killed by a single integer $N$, 
so that 
$$
N(A/(A\cap R))=0.
$$ 
As a consequence, 
$$ 
(A\cap R)[N^{-1}]=A[N^{-1}]. 
$$

Now, let $z\in A$, and let $p$ be a prime. The element 
$$
z\in A[N^{-1}]=(A\cap R)[N^{-1}]
$$ 
can be written $z=\frac{a}{N^k}$, 
where $a\in A\cap R$ and $k\in\N$. Using Fermat's Little Theorem, we find that
$$
(N^k)^p = N^k + pn \qquad\text{and}\qquad a^p = a + pb
$$
for some $n\in \Z$ and $b\in R$. Observe that $pb$ belongs to $A\cap R$, hence to $A[N^{-1}]$, 
so that $b\in A$, as long as $p$ does not divide $N$. We then have
\begin{multline*}
z^p -z = \frac{a^p}{N^{kp}}-\frac{a}{N^k} = \frac{a+pb}{N^k+pn}-\frac{a}{N^k} \\
= \frac{(a+pb)N^k - a(N^k+pn)}{(N^k+pn)N^k} = p\frac{N^kb - na}{(N^k+pn)N^k} = p\frac{N^kb - na}{N^{(p+1)k}}, 
\end{multline*}
and hence $$pu=z^p-z\in A$$ for some $u\in A[N^{-1}]$, assuming $p\nmid N$. But then, in fact, $u\in A$. 

Consequently, for all $z\in A$ and all sufficiently large primes $p$, 
the relation $z^p-z\in pA$ holds, so that $z^p=z$ in $A/pA$. 
Being reduced and artinian (hence semi-simple), $A/pA$ may be written as a product of fields, 
and, because of the equation $z^p=z$, these fields must all 
equal $\Zp$, which means all sufficiently large primes split completely in $A$. 
It will then be a consequence of Chebotarev's Density 
Theorem%
\footnote{(A special case of) Chebotarev's Density Theorem states the following: 
The density of the primes that split completely in a number field 
$K$ equals $\frac{1}{|\mathrm{Gal}(K/\Q)|}$. In our case, this set has density 1.}
that $\Q\otimes R=\Q$. 
Since the assumption is that $R$ is finitely generated as an abelian group, we infer that $R=\Z$.

\emph{Case A2: The ring of fractions is a product of fields.} 
Supposing now the fraction ring of $R$ to be a finite product $\prod K_j$ of fields, the 
projections $R_j$ of $R$ on the factors $K_j$ will each be binomial. 
Hence $R\subseteq \prod R_j$, each $R_j$ being isomorphic to $\Z$, 
according to the above argument.  
But $\Z$ possesses no non-trivial binomial ideals, so, by Goursat's Lemma, 
$R$ must equal a product of copies of $\Z$. 

\emph{Case B: $R$ is infinitely generated as an abelian group.}
Finally, we drop the assumption that $R$ be finitely generated as a group, 
and assume it finitely generated as a ring only. 
The relation $p\mid r^p-r$ ensures that $R/pR$ will be a finitely generated 
torsion group for each prime $p$. It will then have Krull dimension $0$, 
and it follows from Theorem \ref{S: Dimension} that $\dim \Q\otimes R = 0$,
so that $\Q\otimes R$ is a finite-dimensional vector space over $\Q$. 
Only finitely many denominators can be employed in a basis, 
so there exists an integer $M$ for which $R[M^{-1}]$ is finitely generated over $\Z[M^{-1}]$. 

We can now more or less repeat the previous argument. The ring 
$R[M^{-1}]$ will be reduced, and, as before,  
$\Q\otimes R[M^{-1}]=\Q\otimes R$ will be finite-dimensional, hence a product of fields, and we may 
reduce to the case when it is actually a field.

Letting $A$ denote the algebraic integers in  $\Q\otimes R[M^{-1}]$, the factor group $A/(A\cap R)[M^{-1}]$ will be 
finitely generated and torsion, hence killed by some integer,
so that, again, we are led to $(A\cap R)[N^{-1}]=A[N^{-1}]$. 
As before, we may draw the conclusion that $\Q\otimes R=\Q$, and 
consequently that $R=\Z[N^{-1}]$. This concludes the proof of existence. 

\emph{Uniqueness.} 
The ring 
$$
S=\Z[m_1^{-1}]\times \cdots\times \Z[m_k^{-1}],
$$ 
for square-free integers $m_1,\dots,m_k$, is characterised, 
among rings of this same type, by the following two properties. 
\begin{enumerate}
\item There exist $k$ elements $e_1,\dots,e_k\in S$ such that: 
\begin{enumerate}
\item The set $\{e_1,\dots,e_k\}$ is a basis for $\Q\otimes S$.
\item $e_ie_j=\delta_{ij}e_i$ (Kronecker delta).
\end{enumerate}

\item Any such basis  
may be renumbered $\{e_1,\dots,e_k\}$ so that $e_j$ be divisible (in~$S$) by a square-free number
$n$ if and only if $n\mid m_j$.  
\end{enumerate}
The first property shows that the number $k$ is uniquely determined, and the second that different values for $m_i$ 
yield non-isomorphic rings. 
\end{proof}

An alternative proof, in the case that $R$ is finitely generated as an abelian group, runs as follows. 
As an abelian group, $R\cong\Z^n$ since it is torsion-free. Then $R$ embeds as a subring of matrices
$$
R\to \End \Z^n = \Z^{n\times n}, \qquad a\mapsto \lambda_a,
$$
where $\lambda_a$ denotes left multiplication by $a$. Since $R$ is reduced (the first lines of the proof above), 
the minimal polynomial of $a\in R$ has distinct irreducible factors; hence $a$ is diagonalisable over $\C$. 
Because all elements of $R$ commute, they are simultaneously diagonalisable, and so $R$ may be considered 
a subring of $\C^n$. The projection $S$ of $R$ on any of the $n$ factors is a binomial subring of $\C$, 
finitely generated as an abelian group. Now, if $S$ contains a non-integer, it will contain a non-integer $z$ 
with $\Re z>-1$. For such a $z$, we have 
$$
\Abs{\frac{z-n+1}{n}} = \sqrt{1 - \frac{2(\Re z+1)}{n}+ O\left(\frac{1}{n^2}\right) } 
= 1 - \frac{\Re z+1}{n} + O\left(\frac{1}{n^2}\right),
$$
implying that the sequence $\binom{z}{n}\in S$ tends to $0$, without ever reaching it, as $n\to\infty$. 
(We are using here the fact that $\prod_{n=1}^\infty \left(1+\frac{c}{n}\right)=0$ for $c<0$.)
But then clearly $S$ cannot be 
finitely generated as an abelian group. Consequently, $S$ contains only integers, and is either $\Z$ or $0$. 
Again by Goursat's Lemma, $R$ must be a direct product of copies of $\Z$.

\section{Modules}			\label{A: Modules}

By aid of the Structure Theorem, we may classify the 
finitely generated modules over a finitely generated binomial ring 
$$ 
R\cong \Z[m_1^{-1}]\times \cdots\times \Z[m_k^{-1}]
$$ 
as in the next example. 

\begin{example}				\label{X: Module}
It will clearly be sufficient to 
describe the finitely generated modules over $\Z[m^{-1}]$. This is a principal ideal domain, 
whose finitely generated modules admit a decomposition (unique up to reordering) into direct summands
$\Z[m^{-1}]/(q)$,
where $(q)$ is a primary ideal. Such an ideal is either $0$ or a power of a prime, $q=p^n$, in which case 
$$
\Z[m^{-1}]/(p^n) \cong
\begin{cases}
\Z/p^n & \text{if $p\nmid m$;} \\
0 & \text{if $p\mid m$.} 
\end{cases}
$$
Hence a finitely generated module over $\Z[m^{-1}]$ is isomorphic to a unique module of the form 
$$
\Z[m^{-1}]^r \oplus \Z/p_1^{n_1} \oplus \cdots \oplus \Z/p_l^{n_l},
$$
where neither of the primes $p_1,\dots,p_l$ divides $m$.
\end{example}

Without too much labour, one may arrive at a slightly more general proposition, classifying 
torsion-free modules $M$ over a binomial ring $R$, with the property that 
 $M[m^{-1}]$ is finitely generated over $\Z[m^{-1}]$ for some integer $m$.
Note that, this time, we do not assume $R$ to be finitely generated. 

Our first observation is that $\End M[m^{-1}]$ will also be a finitely generated module over $\Z[m^{-1}]$
and thus noetherian, so that 
the submodule $R[m^{-1}]$ is finitely generated as well.
\emph{A fortiori,} $R[m^{-1}]$ is finitely generated as a ring, 
and, using the Structure Theorem, we deduce its form to be
$$ 
R[m^{-1}]\cong \Z[n_1^{-1}]\times \cdots\times \Z[n_k^{-1}]. 
$$ 
For this to be a $\Z[m^{-1}]$-module, we must have $m\mid n_i$, and for it to be finitely generated as such, 
we must have $n_i=m$; hence $R[m^{-1}]\cong \Z[m^{-1}]^k$.

Next, the projection of $R$ on the $i$'th factor must be of the form $\Z[m_i^{-1}]$, with $m_i$ a factor of $m$.
Therefore, 
$$ 
R \subseteq \Z[m_1^{-1}]\times \cdots\times \Z[m_k^{-1}],
$$
the projection on each factor now being surjective. 
Since each factor $\Z[m_i^{-1}]$ is devoid of non-trivial binomial ideals, 
Goursat's Lemma applies once more to show that 
$R$ itself is, in fact, a product of such rings (possibly containing fewer than $k$ direct factors).

We have thus proved: 

\begin{theorem}
Let $M$ be a torsion-free module over the binomial ring $R$, with the property that 
$M[m^{-1}]$ is finitely generated over $\Z[m^{-1}]$ for some integer $m$.
Then 
$$
R \cong \Z[m_1^{-1}]\times \cdots\times \Z[m_k^{-1}]
$$
for some integers $m_i\mid m$, and $M$ admits of an explicit description as in Example \ref{X: Module}.
\end{theorem}

\end{document}